\documentclass[reqno,12pt]{amsart}

\usepackage{enumerate}
\usepackage[margin=1in]{geometry}
\usepackage{ifpdf}
\usepackage{amsmath}
\usepackage{amsfonts}
\usepackage{amssymb}
\usepackage{amsthm}
\usepackage[hyperfootnotes=false]{hyperref}
\usepackage{setspace}
\usepackage{amsrefs}
\usepackage{nicefrac}
\usepackage{graphicx}

\newcommand{\ignore}[1]{}

\newcommand{\C}{{\mathbb{C}}}
\newcommand{\R}{{\mathbb{R}}}

\newcommand{\sS}{{\mathcal{S}}}

\newtheorem{thm}{Theorem}

\newtheorem{prop}[thm]{Proposition}

\theoremstyle{definition}

\theoremstyle{remark}

\author{Ji\v{r}\'{\i} Lebl}
\thanks{The first author was in part supported by NSF grant DMS-1362337 and
Oklahoma State University's DIG and ASR grants.}
\address{Department of Mathematics, Oklahoma State University,
Stillwater, OK 74078, USA}
\email{lebl@math.okstate.edu}
\date{April 14, 2021}

\ifpdf
\hypersetup{
  pdftitle={An example of a compact non-$\C$-analytic real subvariety of $\R^3$},
  pdfauthor={Jiri Lebl},
  pdfsubject={Real analytic varieties},
  pdfkeywords={Real analytic varieties},
}
\fi

\title%
{An example of a compact non-$\C$-analytic real subvariety of $\R^3$}

\begin{document}


\begin{abstract}
The purpose of this short expository note is to provide an example exhibiting some of the pathological properties of real-analytic subvarieties, where the pathology can be visualized, and the proofs use only elementary properties of analytic functions.  We construct a compact irreducible real-analytic subvariety $S$ of ${\mathbb R}^3$ of pure dimension two such that 1) the only a real-analytic function is defined in a neighbourhood of $S$ and vanishing on $S$ is the zero function, 2) the singular set of $S$ is not a subvariety of $S$, nor is it contained in any one-dimensional subvariety of $S$, 3) the variety $S$ contains a proper subvariety of dimension two.  The example shows how a badly behaved part of a subvariety can be hidden via a second well-behaved component to create a subvariety of a larger set.  The pathology is visualized using several figures.  Examples of these phenomena are known since the time of Cartan, but hard to find in the English language literature.
\end{abstract}

\maketitle




A general problem in geometry is to study the set of solutions to a certain
equation.  Classically one studies polynomial equations, and then one
is led to extend this study to analytic equations.  Over the
complex field, these solution sets are well-behaved due to the
algebraic completeness of the complex numbers.  Over the real numbers,
a rather peculiar thing happens.  Even though real-analytic functions are
(locally) restrictions of complex-analytic functions, their solutions
exhibit pathologies not present in the complex case.

A closed set $S \subset \R^n$ is a \emph{real-analytic subvariety} if 
for every $p \in S$, there exist real-analytic functions
$\rho_1,\ldots,\rho_k$ defined in a neighbourhood $U$ of $p$, such that $S
\cap U$ is equal to the set where all $\rho_1,\ldots,\rho_k$ vanish.
A complex subvariety is the same notion with
$\R^n$ replaced by $\C^n$ and
real-analytic replaced by holomorphic (complex analytic).
See \cites{GMT:topics,Whitney:book} for more information.

Given a real-analytic function $\rho(x_1,\ldots,x_n)$ defined near $p$, 
consider $\R^n \subset \C^n$, and then consider $x_1,\ldots,x_n$
as complex variables.  The power series still converges in some small
neighbourhood of $p \in \C^n$.  This process is called
\emph{complexification}.
Starting with a
real-analytic subvariety near a point $p$, complexify the
defining equations near $p$
to obtain a complex subvariety of a neighbourhood of $p \in \C^n$ whose
trace on $\R^n$ is the real subvariety.
Such a complex subvariety may not exist globally, but only
locally.  A subvariety in $\R^n$ which is the real trace of a 
complex variety is usually called $\C$-analytic,
and the purpose here is to provide a non-$\C$-analytic example
subvariety.

The basic properties of real-subvarieties including examples of the
pathologies explored in this note have been known at least since the work of
Cartan~\cite{Cartan}, see also \cites{BruhatCartan,WhitneyBruhat}.
There exist examples in the literature, but they may be hard to find,
and not well-known outside the field,
especially since many appear in non-English literature.
There are only few accessible books on the subject,
e.g.~\cite{GMT:topics,narasimhan}, which seem to be missing
in many university libraries.
This note is hoped to improve the situation.
Furthermore, just in the recent years the author can count at least two
preprints that included claims which are refuted by the existence of
these examples, so it appears such a note is necessary.
One other motivation for this note is to present an example subvariety
where the pathology can be easily visualised.
In particular, we can visualize how a badly behaved part of a subvariety
can be hidden to
extend the subvariety to a larger set via a seemingly unrelated other
subvariety with a very different geometry.  

A real-analytic subvariety
$X \subset \R^n$ is \emph{irreducible} if whenever we write 
$X = X_1 \cup X_2$ for two subvarieties $X_1$ and $X_2$ of $\R^3$, then
either $X_1 = X$ or $X_2 = X$.  This notion is subtle.  We will construct
a set that is a union of two subvarieties, one of which
is not a subvariety of $\R^3$ but of a strictly smaller domain, and 
the union is irreducible as a subvariety of $\R^3$.

\begin{figure}[ht]
\hfill
\includegraphics[width=2in]{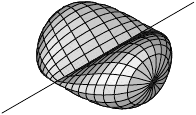}
\hfill
\includegraphics[width=2.5in]{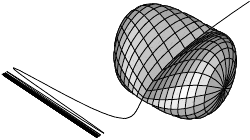}
\hfill{}
\caption{The sets $S_1$ (left) and $S_2$ (right).\label{fig:bubble12}}
\end{figure}

Let us proceed in steps.
Start with the sphere $z^2=1-x^2-y^2$, thinking
of $z^2$ as a ``graph''.  Pinch the sphere along the $y$-axis by
multiplying by $x^2$ to obtain the subvariety $S_1$ given by
\begin{equation} \label{eq:S1eq}
z^2=(1-x^2-y^2)x^2 .
\end{equation}
The picture is the left hand side of Figure~\ref{fig:bubble12}.

The subvariety
$S_1$ is irreducible, and it contains
the $y$-axis as a subvariety.  This subvariety is already somewhat
pathological.  It has components of different dimensions, and the regular
points of dimension two are not topologically connected.
By \emph{regular} points we mean points near which the subvariety is a
real-analytic submanifold.
It can be shown that the regular points are dense in the subvariety,
and by analytic continuation, if the set of regular points is connected,
then the subvariety must be irreducible.  The converse is true
(but much harder to prove) for complex subvarieties, but as we see,
it is not true for real subvarieties.

We restrict our attention to the set where $-2 < y < 2$.  We 
change coordinates
keeping $y$ and $z$ fixed, but sending
$x$ to $x+\sin\bigl(\frac{1}{y+2}\bigr)$.  The equation becomes
\begin{equation} \label{eq:S2eq}
z^2=\left(1-{\left(x-\sin\left(\frac{1}{y+2}\right)\right)}^2-y^2\right)
{\left(x-\sin\left(\frac{1}{y+2}\right)\right)}^2.
\end{equation}
Call this set $S_2$, restricted to $-2 < y < 2$ to make it bounded.
What we have done is made the $y$ axis appear like the graph of 
$\sin\bigl(\frac{1}{y+2}\bigr)$.  The picture of the result is
the right hand side of Figure~\ref{fig:bubble12}.
The set $S_2$ is a real-analytic subvariety of the set $\{ -2 < y < 2 \}$,
and it becomes badly behaved as we approach $y=-2$.

\begin{prop}
Let 
$\Omega \subset \R^3$ be a connected neighbourhood of the closure $\overline{S_2}$
and $r \colon \Omega \to \R$ a real-analytic function vanishing on
$S_2$.
Then $r$ is identically zero.
\end{prop}

It is easy to prove that $r$ must vanish on the $xy$-plane
by staring at the picture and recalling basic
properties of real-analytic functions.  To prove that
$r$ vanishes everywhere, we need to complexify $r$.
The result is not simply because the 1-dimensional component wiggles
around, it is because $r$ must vanish on the ``complexification'' of the
2-dimensional part,
which gets ``dragged along'' the 1-dimensional component.

\begin{proof}
We consider $\R^3 \subset \C^3$.  Then $S_2$ is also a subset of $\C^3$.
Treating $(x,y,z)$ as complex variables,
let us call $\sS$ the complex subvariety of $\{ y \not= -2 \}$ set defined by
\eqref{eq:S2eq}.  
The subvariety $\sS$ is locally irreducible at $(0,-1,0)$.
Indeed, think of $z^2$ as a graph over $(x,y)$, and so there can at most be ``two sheets'' in $\sS$, for the two
different square roots of the right hand side of \eqref{eq:S2eq}.
The set of points where the right hand side of \eqref{eq:S2eq}
does not vanish correspond to regular points of $\sS$.
The zero set of the right hand side (considering $(x,y) \in \C^2$)
is a complex subvariety itself and thus its complement is connected.
As the set of regular points is connected,
the subvariety is irreducible as we mentioned above.

Complexify $r$ to obtain a holomorphic
function $\widetilde{r}$ of a neighbourhood $U$ of $\overline{S_2}$ in $\C^3$.
As $\widetilde{r}$ vanishes on $\overline{S_2}$, then by irreducibility of $\sS$
at $(0,-1,0)$, there exists a neighbourhood $W$ of $(0,-1,0)$ in $\C^3$,
such that $\widetilde{r}$ vanishes on $W \cap \sS$.

Fix $z = ia$ for a small real $a$.  Let $X_a$ be the set defined by
\begin{equation} \label{eq:Xaeq}
-a^2=\left(1-{\left(x-\sin\left(\frac{1}{y+2}\right)\right)}^2-y^2\right)
{\left(x-\sin\left(\frac{1}{y+2}\right)\right)}^2.
\end{equation}
for real $x$ and $y$ with $-2 < y < -1$.  The set $X_a$ is a connected
smooth real-analytic curve, which is a subset of $\sS$.
If $a$ is small enough, 
$X_a \subset \sS \cap U$ and $X_a \cap W$ is nonempty.
The function $\widetilde{r}$ then
vanishes on $X_a \cap W$, an open set of $X_a$, and hence on all of $X_a$.

Next fix a small real $x$.  The equation \eqref{eq:Xaeq} is true for an
infinite sequence of $y$ approaching $y=-2$ from above.  Therefore, the
holomorphic function of one complex variable $y \mapsto
\widetilde{r}(x,y,ia)$ defined in a neighbourhood of the origin
vanishes identically.  As this was true for all small enough $x$ and $a$,
it is true for small enough complex $x$ and $a$,
and $\widetilde{r}$ vanishes in a neighbourhood of
$(0,-2,0)$ in $\C^3$.  By analytic continuation $\widetilde{r}$ is
identically zero.
\end{proof}

\begin{figure}[ht]
\hfill
\includegraphics[width=2.3in]{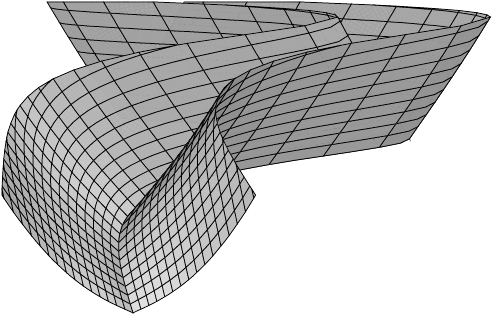}
\hfill{}
\caption{The trace of the complex variety containing $S_2$ in $(x,y,a)$-space for $a \geq 0$
and $y < -1$,
approaching $y=-2$.\label{fig:valley}}
\end{figure}

To visualize how bad the complex variety containing $S_2$ is as we approach $y=-2$, consider
the set in the space $(x,y,a) \in \R^3$ given by \eqref{eq:Xaeq}.  Looking at the set
where $a \geq 0$,
we have a ``valley'' whose bottom is the graph 
$x = \sin\bigl(\frac{1}{y+2}\bigr)$ with increasingly steep sides.  See
Figure~\ref{fig:valley}.

\begin{figure}[ht]
\hfill
\includegraphics[width=2.5in]{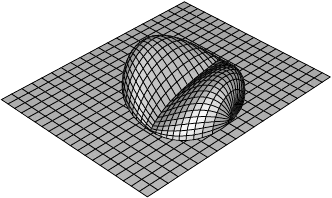}
\hfill
\includegraphics[width=2.0in]{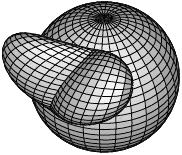}
\hfill{}
\caption{The sets $S_3$ (left) and $S_4$ (right).\label{fig:bubble34}}
\end{figure}

Let us ``hide'' the wild behavior near $y=-2$ and
construct the purely 2-dimensional subvariety $S_3$ via
\begin{equation}
S_3 = S_2 \cup \{ z = 0 \} .
\end{equation}
The picture is the left hand side of Figure~\ref{fig:bubble34}.
Suppose $\Omega \subset \R^3$ is a connected
neighbourhood of $S_3$ and $r \colon \Omega \to \R$
a real-analytic function such that $r=0$ on $S_3$.  The set $\Omega$ is also a neighbourhood of
$\overline{S_2}$ and $r=0$ on $\overline{S_2}$.
By the proposition, $r \equiv 0$.
In the terminology of real-analytic varieties, $S_3$ is not $\C$-analytic.

The singular set of $S_3$ is the set
\begin{equation}
\begin{aligned}
z & = 0 ,  \\
-1 & \leq y \leq 1 , \\
0 & =\left(1-{\left(x-\sin\left(\frac{1}{y+2}\right)\right)}^2-y^2\right)
{\left(x-\sin\left(\frac{1}{y+2}\right)\right)}^2 .
\end{aligned}
\end{equation}
This singular set clearly is not a subvariety.  In particular, it contains the set
\begin{equation}
I = \left\{ (x,y,z) :
x=\sin\left(\frac{1}{y+2}\right) \quad \text{and} \quad -1 \leq y \leq 1
\right\},
\end{equation}
and $I$ cannot be contained in any subvariety of $S_3$
of dimension 1.
Any such subvariety would have to contain the entire set
$x=\sin\bigl(\frac{1}{y+2}\bigr)$ for all $y > -2$, and it cannot
possibly
be a subvariety at
points where $y = -2$.

The subvariety $S_3$ is irreducible.  It contains a proper subvariety of 
dimension 2, namely the $xy$-plane.  Any subvariety $S'$ that contains any open
set of the regular points must contain $I$.  Indeed, if $S'$ contains an open
set of the $xy$-plane it must contain whole $xy$-plane.  If $S'$ contains
an open set of one of the smooth submanifolds outside of the $xy$-plane then
$I$ is in the closure of this submanifold and hence in $S'$.  Any subvariety
that contains $I$ must contain the entire $xy$-plane.

We have demonstrated a subvariety with all the required properties but not a
compact one.  We make the subvariety compact by mapping the plane
onto the sphere using spherical coordinates.  For the picture on
the right hand side of Figure~\ref{fig:bubble34} we used the map
\begin{equation}
(x,y,z) \mapsto \Bigl(
(z+1)\sin(1+y/2)\cos(x),
(z+1)\sin(1+y/2)\sin(x),
(z+1)\cos(1+y/2) \Bigr) .
\end{equation}
For $-\pi < x < \pi$, $0 < (1+y/2) < \pi$, and $z+1 > 0$, the
mapping is a real-analytic diffeomorphism and we obtain the compact
subvariety $S_4$ by taking the closure, which will fill in the missing
meridian on the far side of the sphere.

This subvariety clearly has all the properties mentioned in the abstract;
it is compact and inherits the rest of the properties from $S_3$.

The construction is easy to modify to show further strange behaviors.
For example, if we start with $S_3$, but rescale the $z$
variable we obtain another irreducible
subvariety that shares with $S_3$ a 2-dimensional component as a proper
subvariety.


\def\MR#1{\relax\ifhmode\unskip\spacefactor3000 \space\fi%
  \href{http://www.ams.org/mathscinet-getitem?mr=#1}{MR#1}}

\end{document}